\documentclass[10pt]{amsart}

\usepackage{amssymb,amsmath,amsthm,amsfonts}

\theoremstyle{plain}
\newtheorem{Theorem}{Theorem}

\newtheorem{Conjecture}[Theorem]{Conjecture}

\theoremstyle{definition}

\begin{document}
\title[Abelianization conjectures]{Abelianization conjectures for some arithmetic square complex groups}
\author{Diego Rattaggi}
\thanks{Supported by the Swiss National Science Foundation, No.\ PP002--68627}
\address{Universit\'e de Gen\`eve,
Section de math\'ematiques,
2--4 rue du Li\`evre, CP 64, 
CH--1211 Gen\`eve 4, Switzerland}
\email{rattaggi@math.unige.ch}
\date{\today}
\begin{abstract}
We extend a conjecture of Kimberley-Robertson on the abelianizations of certain
square complex groups.
\end{abstract}
\maketitle
\section{Introduction} \label{Intro}
Throughout this paper, let $p$, $l$ be any pair of distinct odd prime numbers,
\[
r_{p,l} := \gcd \left( \frac{p-1}{4}, \frac{l-1}{4}, 6 \right) \in \{ 1,2,3,6 \} \, ,
\]
and $q \in \{p,l\}$.
We first recall the definition of the group $\Gamma_{p,l}$ from
\cite{Rattaggi,Rattaggi2,Rattaggi3,RaRo}.
Let $\mathbb{Q}_q$ be the field of $q$-adic numbers.
We fix elements 
$c_p, d_p \in \mathbb{Q}_p$ and $c_l, d_l \in \mathbb{Q}_l$
such that 
\[
c_p^2 + d_p^2 + 1 = 0 \in \mathbb{Q}_p \,
\text{ and } \,
c_l^2 + d_l^2 + 1 = 0 \in \mathbb{Q}_l.
\]
Note that we can take $d_q = 0$, if $q \equiv 1 \pmod{4}$.

Let $\mathbb{H}(\mathbb{Q})^{\ast}$ be the multiplicative group of
invertible rational Hamilton quaternions, i.e.\ the set
\[
\{x_0 + x_1 i + x_2 j + x_3 k :  x_0, x_1, x_2, x_3 \in \mathbb{Q} \} \setminus \{ 0 \}
\]
equipped with the multiplication induced by the rules $i^2 = j^2 = k^2 = -1$ and $ij = k = -ji$.
If $x = x_0 + x_1 i + x_2 j + x_3 k$, we define as usual 
the conjugate $\overline{x} := x_0 - x_1 i - x_2 j - x_3 k$,
and the norm $|x|^2 := x \overline{x} = x_0^2 + x_1^2 + x_2^2 + x_3^2$.

Let $\psi_q$ be the homomorphism of groups
$\mathbb{H}(\mathbb{Q})^{\ast} \to \mathrm{PGL}_2(\mathbb{Q}_q)$
defined by 
\[
\psi_q(x_0 + x_1 i + x_2 j + x_3 k) := \left[
\begin{pmatrix}
x_0 + x_1 c_q + x_3 d_q & -x_1 d_q + x_2 + x_3 c_q \\
-x_1 d_q - x_2 + x_3 c_q  & x_0 - x_1 c_q - x_3 d_q \\
\end{pmatrix}\right]
\]
and let the homomorphism
\[
\psi_{p,l} : \mathbb{H}(\mathbb{Q})^{\ast} \to 
\mathrm{PGL}_2(\mathbb{Q}_p) \times \mathrm{PGL}_2(\mathbb{Q}_l)
\]
be given by $\psi_{p,l}(x) := (\psi_p(x), \psi_l(x))$.
Observe that it satisfies $\psi_{p,l}(-x) = \psi_{p,l}(x)$ and 
$\psi_{p,l}(x)^{-1} = \psi_{p,l}(\overline{x})$.

Let $\mathbb{H}(\mathbb{Z})$ be the set of integer Hamilton quaternions and
$X_q$ the subset of quaternions
\begin{align*}
X_q := \{x = \, & x_0  + x_1 i + x_2 j + x_3 k \in \mathbb{H}(\mathbb{Z}) \,; 
\quad |x|^2 = q \, ; \notag \\
& x_0 \text{ odd}, \text{ if } q \equiv 1 \!\!\!\! \pmod 4 \, ; \;
  x_1 \text{ even}, \text{ if } q \equiv 3 \!\!\!\! \pmod 4 \} \notag
\end{align*}
of cardinality $2(q+1)$.

Finally, let $Q_{p,l}$ be the subgroup of $\mathbb{H}(\mathbb{Q})^{\ast}$
generated by $(X_p \cup X_l) \subset \mathbb{H}(\mathbb{Z})$ and let 
$\Gamma_{p,l} < \mathrm{PGL}_2(\mathbb{Q}_p) \times \mathrm{PGL}_2(\mathbb{Q}_l)$ 
be its image $\psi_{p,l}(Q_{p,l})$, which is a finitely presented, torsion-free linear group. 

The starting point for this work was the following conjecture of Kimberley and Robertson 
for the abelianization $\Gamma_{p,l}^{ab} := \Gamma_{p,l} / [\Gamma_{p,l}, \Gamma_{p,l}]$
of the group $\Gamma_{p,l}$ in the case $p,l \equiv 1 \pmod{4}$.
We use the notation $\mathbb{Z}_n := \mathbb{Z} / n \mathbb{Z}$ and
$\mathbb{Z}_n^m := \mathbb{Z} / n \mathbb{Z} \times \ldots \times \mathbb{Z} / n \mathbb{Z}$
($m$ times).

\begin{Conjecture} \label{KimbRob}
(Kimberley-Robertson \cite[Section~6]{KR}) 
If $p,l \equiv 1 \pmod{4}$, then
\[
\Gamma_{p,l}^{ab} \cong
\begin{cases}
\mathbb{Z}_2 \times \mathbb{Z}_4^3 \, , & \text{if } \, r_{p,l} = 1 \\
\mathbb{Z}_2^3 \times \mathbb{Z}_8^2 \, , & \text{if } \, r_{p,l} = 2 \\
\mathbb{Z}_2 \times \mathbb{Z}_3 \times \mathbb{Z}_4^3 \, , & \text{if } \, r_{p,l} = 3 \\
\mathbb{Z}_2^3 \times \mathbb{Z}_3 \times \mathbb{Z}_8^2 \, , & \text{if } \, r_{p,l} = 6 \, .
\end{cases}
\]
\end{Conjecture}
In Section~\ref{CaseOne}, we will give an equivalent formulation of this conjecture 
and a new conjecture relating the abelianization of $\Gamma_{p,l}$ to the number 
$t_{p,l}$ of certain pairs of commuting quaternions, defined as
\[
t_{p,l} := | \{ (x,y) \in Y_p \times Y_l \, : \, xy = yx \} | \, ,
\]
where $Y_q$ is any subset of $X_q$ of cardinality $(q+1)/2$ such that $x \in Y_q$ implies
$x_0 > 0$ and $\overline{x} \notin Y_q$. 
Note that the definition of $t_{p,l}$ does not depend on the choice of elements in $Y_p$ and $Y_l$,
and that $\psi_{p,l}(Y_p \cup Y_l)$ is a generating set of $\Gamma_{p,l}$
of cardinality $(p+1)/2 + (l+1)/2$.

The case $p, l \equiv 3 \pmod{4}$ is treated in Section~\ref{CaseTwo}
and the remaining (mixed) case in Section~\ref{CaseThree}.
The final section is devoted to conjectures on the abelianizations
of some subgroups of $\Gamma_{p,l}$.

The Conjectures~\ref{KimbRobEquiv},
\ref{Conjcpl},
\ref{ConjcplKimbRob},
\ref{Conj334},
\ref{Conj314},
\ref{ConjComm114},
\ref{ConjComm334},
\ref{ConjComm314} and
\ref{ConjectureLambda} have been
stated in the authors Ph.D.~thesis (\cite[Chapter~3]{Rattaggi}).
We have checked Conjectures~\ref{KimbRobEquiv},
\ref{ConjcplKimbRob},
\ref{Conj334},
\ref{Conj314},
\ref{ConjComm114},
\ref{ConjComm334},
\ref{ConjComm314}
and
\ref{ConjectureLambda}
for more than $100$ different pairs $(p,l)$ which are explicitly listed in
\cite[Table~3.13]{Rattaggi}. 

\section{The case $p,l \equiv 1 \pmod{4}$} \label{CaseOne}
In this section, we restrict to the ``classical'' case $p,l \equiv 1 \pmod{4}$. The following conjecture
is equivalent to Conjecture~\ref{KimbRob}.
\begin{Conjecture} \label{KimbRobEquiv}
Let $p, l \equiv 1 \pmod{4}$.\\
If $p, l \equiv 1 \pmod{8}$, then
\[
\Gamma_{p,l}^{ab} \cong
\begin{cases}
\mathbb{Z}_2^3 \times \mathbb{Z}_3 \times \mathbb{Z}_8^2 \, , & \text{if }  \, p, l \equiv 1 \!\!\! \pmod{3} \\
\mathbb{Z}_2^3 \times \mathbb{Z}_8^2 \, ,                     & \text{else} \, .
\end{cases}
\]
If $p  \equiv 5 \pmod{8}$ or $l \equiv 5 \pmod{8}$, then
\[
\Gamma_{p,l}^{ab} \cong
\begin{cases}
\mathbb{Z}_2 \times \mathbb{Z}_3 \times \mathbb{Z}_4^3 \, , & \text{if }  \, p, l \equiv 1 \!\!\! \pmod{3} \\
\mathbb{Z}_2 \times \mathbb{Z}_4^3 \, ,                     & \text{else} \, .
\end{cases}
\]
\end{Conjecture}

\begin{proof}[Proof of the equivalence of Conjecture~\ref{KimbRob} and Conjecture~\ref{KimbRobEquiv}]
If $r_{p,l} = 6$, then $(p-1)/4 = 6s$ and $(l-1)/4 = 6t$ for some $s, t \in \mathbb{N}$, i.e.\
$p = 24s + 1$ and $l = 24t +1$. It follows $p, l \equiv 1 \pmod{8}$ and $p, l \equiv 1 \pmod{3}$.

If $r_{p,l} = 3$, then $(p-1)/4 = 3s$ and $(l-1)/4 = 3t$, where $s$ or $t$ is odd 
(otherwise~$r_{p,l}$ would be $6$). 
Consequently, we have $p = 12s + 1$ and $l = 12t +1$,
in particular $p, l \equiv 1 \pmod{3}$. If $s$ is odd, then $p \equiv 5 \pmod{8}$. 
If $t$ is odd, then $l \equiv 5 \pmod{8}$.

If $r_{p,l} = 2$, then $(p-1)/4 = 2s$ and $(l-1)/4 = 2t$, i.e.\ $p = 8s + 1$ and $l = 8t +1$, 
hence $p, l \equiv 1 \pmod{8}$.
Moreover, $s \not\equiv 0 \pmod{3}$ or $t \not\equiv 0 \pmod{3}$ (otherwise~$r_{p,l}$ would be $6$). 
In the first case, we have 
$p \not\equiv 1 \pmod{3}$, in the second case $l \not\equiv 1 \pmod{3}$.

If $r_{p,l} = 1$, then $(p-1)/4 = 2s - 1$ or  $(l-1)/4 = 2t - 1$ (otherwise $r_{p,l}$ would be even), 
hence 
$p = 8s - 3$ or $l = 8t - 3$, i.e.\ $p \equiv 5 \pmod{8}$ or $l \equiv 5 \pmod{8}$.
Moreover: $(p-1)/4 = 3s + 1$ or $(p-1)/4 = 3s + 2$ or $(l-1)/4 = 3t + 1$  or $(l-1)/4 = 3s + 2$ 
for some $s,t \in \mathbb{N}_0$
(otherwise $r_{p,l}$ would be a multiple of $3$), hence $p = 12s + 5$ or $p = 12s + 9$ or 
$l = 12t + 5$ or $l = 12t + 9$,
in particular $p \not\equiv 1 \pmod{3}$ or $l \not\equiv 1 \pmod{3}$.
\end{proof}
The equivalence of the two conjectures above is also expressed in Table~\ref{Table24}.
\begin{table}[ht]
\[
\begin{tabular}{r|rccccl}
$r_{p,l}$    & $l \equiv 1$ & $5$ & $9$ & $13$ & $17$ & $21 \!\!\! \pmod{24}$ \notag \\ \hline
$p \equiv 1$ & $6$          & $1$ & $2$ & $3$  & $2$  & $1$ \notag \\ 
$5$          & $1$          & $1$ & $1$ & $1$  & $1$  & $1$ \notag \\
$9$          & $2$          & $1$ & $2$ & $1$  & $2$  & $1$ \notag \\
$13$         & $3$          & $1$ & $1$ & $3$  & $1$  & $1$ \notag \\
$17$         & $2$          & $1$ & $2$ & $1$  & $2$  & $1$ \notag \\
$21$         & $1$          & $1$ & $1$ & $1$  & $1$  & $1$ \notag 
\end{tabular}
\]   
\caption{$r_{p,l}$ for $p$, $l$ taken modulo $24$} 
\label{Table24}
\end{table}

The structure of $\Gamma_{p,l}^{ab}$ also seems to depend only on the number $t_{p,l}$
defined in Section~\ref{Intro}.
Observe that 
\[
3 \leq t_{p,l} \leq \mathrm{min} \left\{ \frac{p+1}{2}, \frac{l+1}{2} \right\} \, ,
\]
if $p,l \equiv 1 \pmod{4}$.
\begin{Conjecture} \label{Conjcpl}
Let $p,l \equiv 1 \pmod{4}$. Then
\[
t_{p,l} \equiv
\begin{cases}
3 \!\!\!\!\pmod{12} \, , & \text{if } \, r_{p,l} = 1 \\
9 \!\!\!\!\pmod{12} \, , & \text{if } \, r_{p,l} = 2 \\
7 \!\!\!\!\pmod{12} \, , & \text{if } \, r_{p,l} = 3 \\
1 \!\!\!\!\pmod{12} \, , & \text{if } \, r_{p,l} = 6 \, .
\end{cases}
\]
\end{Conjecture}
We have checked Conjecture~\ref{Conjcpl} for all pairs of distinct
prime numbers $p,l < 1000$ such that $p,l \equiv 1 \pmod{4}$.
The following values for $t_{p,l}$ appear in this range:
\[
t_{p,l} \in
\begin{cases}
\{ 3,15,27,39,51,63,75,87,99\} \, , & \text{if } \, r_{p,l} = 1 \\
\{ 9,21,33,45,57,69,81,93,105,117,129,153\} \, , & \text{if } \, r_{p,l} = 2 \\
\{ 7,19,31,43,55,67,79,91,103,115,127,151\} \, , & \text{if } \, r_{p,l} = 3 \\
\{ 37,49,61,73,85,97,109,121,133\} \, , & \text{if } \, r_{p,l} = 6 \, .
\end{cases}
\]
See Table~\ref{Tablecplfreq} 
for the frequencies of the values of $t_{p,l}$, where $p,l \equiv 1 \pmod{4}$
are prime numbers such that $p < l < 1000$.
\begin{table}[ht]
\[
\begin{tabular}{c|rrrrrrr||r}
$t_{p,l}$ & \textbf{3}   & \textbf{15}  & \textbf{27} & \textbf{39} & \textbf{51} & \textbf{63} & \textbf{75} \notag \\
$\#$      &1242 &449  &143 & 56 & 34 & 17 &  7 \notag \\ \hline
          & \textbf{87}  & \textbf{99}  &    &     &     &      && \notag \\
          &  5  &  2  &    &     &     &      && 1955 \notag \\ \hline \hline
$t_{p,l}$ & \textbf{9}  & \textbf{21} & \textbf{33} & \textbf{45} & \textbf{57} & \textbf{69} & \textbf{81} & \notag \\
$\#$      &178 &158 & 84 & 57 & 40 & 21 &  8 & \notag \\ \hline
          & \textbf{93} & \textbf{105} & \textbf{117} & \textbf{129} & 141 & \textbf{153}  && \notag \\
          &  9 &  12 &   5 &   2 &     &   1  && 575 \notag \\ \hline \hline
$t_{p,l}$ & \textbf{7}  & \textbf{19} & \textbf{31} & \textbf{43} & \textbf{55} & \textbf{67} & \textbf{79} & \notag \\
$\#$      &236 &130 & 79 & 42 & 18 &  8 & 12 & \notag \\ \hline
          & \textbf{91} & \textbf{103} & \textbf{115} & \textbf{127} & 139 & \textbf{151}  && \notag \\
          &  6 &   1 &   4 &   2 &     &   1  && 539 \notag \\ \hline \hline
$t_{p,l}$ & 1  & 13 & 25 & \textbf{37} & \textbf{49} & \textbf{61} & \textbf{73} & \notag \\
$\#$      &    &    &    & 26 & 15 & 15 & 16 & \notag \\ \hline
          & \textbf{85} &  \textbf{97} & \textbf{109} & \textbf{121} & \textbf{133} &      && \notag \\
          & 7 &   4 &   3 &   2 &   3 &       && 91 \notag \\ \hline \hline
          &    &    &    &    &    &    &        & 3160
\end{tabular}
\]   
\caption[Number $t_{p,l}$ and its frequency]{$t_{p,l}$ and its frequency, $p < l < 1000$} \label{Tablecplfreq}
\end{table}

Combining Conjecture~\ref{Conjcpl} with Conjecture~\ref{KimbRob}, we get another conjecture:
\begin{Conjecture} \label{ConjcplKimbRob}
Let $p,l \equiv 1 \pmod{4}$. Then
\[
\Gamma_{p,l}^{ab} \cong
\begin{cases}
\mathbb{Z}_2 \times \mathbb{Z}_4^3 \, , & \text{if } \, t_{p,l} \equiv 3 \!\!\! \pmod{12} \\
\mathbb{Z}_2^3 \times \mathbb{Z}_8^2 \, , & \text{if } \, t_{p,l} \equiv 9 \!\!\! \pmod{12} \\
\mathbb{Z}_2 \times \mathbb{Z}_3 \times \mathbb{Z}_4^3 \, , & \text{if } \, t_{p,l} \equiv 7 \!\!\! \pmod{12} \\
\mathbb{Z}_2^3 \times \mathbb{Z}_3 \times \mathbb{Z}_8^2 \, , & \text{if } \, t_{p,l} \equiv 1 \!\!\! \pmod{12}  \, .
\end{cases}
\]
\end{Conjecture}

\section{The case $p,l \equiv 3 \pmod{4}$} \label{CaseTwo}
If $p,l \equiv 3 \pmod{4}$, we have a conjecture similar to Conjecture~\ref{KimbRobEquiv}.
\begin{Conjecture} \label{Conj334}
Let $p,l \equiv 3 \pmod{4}$.\\
If $p \! \pmod{8} \, = \, l  \! \pmod{8}$, then
\[
\Gamma_{p,l}^{ab} \cong
\begin{cases}
\mathbb{Z}_2 \times \mathbb{Z}_3 \times \mathbb{Z}_8^2 \, , & 
\text{if  } \, p, l \equiv 1 \!\!\! \pmod{3} \; \text{\phantom{iii}(=: case (B1))}\\
\mathbb{Z}_2 \times \mathbb{Z}_8^2 \, ,                     & 
\text{else} \; \text{\phantom{iiiiiiiiiiiiiiiiiiiiiiiii}(=: case (B2))} \, .
\end{cases}
\]
If $p  \! \pmod{8} \, \ne \, l  \! \pmod{8}$, then
\[
\Gamma_{p,l}^{ab} \cong
\begin{cases}
\mathbb{Z}_2 \times \mathbb{Z}_3 \times \mathbb{Z}_4^2 \, , & 
\text{if } \, p, l \equiv 1 \!\!\! \pmod{3} \; \text{\phantom{iii}(=: case (B3))}\\
\mathbb{Z}_2 \times \mathbb{Z}_4^2 \, ,                     & 
\text{else} \; \text{\phantom{iiiiiiiiiiiiiiiiiiiiiiiii}(=: case (B4))} \, .
\end{cases}
\]
\end{Conjecture}
The four cases (B1)--(B4) defined in the conjecture above can also be 
expressed taking $p$ and $l$ modulo $24$, see Table~\ref{Table33424}.
\begin{table}[ht]
\[
\begin{tabular}{r|rccccl}
             & $l \equiv 3$ & $7$ & $11$ & $15$ & $19$ & $23 \!\!\! \pmod{24}$ \notag \\ \hline
$p \equiv 3$ & (B2)         & (B4)&  (B2)&  (B4)&   (B2)&  (B4)\notag \\ 
$7$          & (B4)         & (B1)&  (B4)&  (B2)&   (B3)&  (B2)\notag \\
$11$         & (B2)         & (B4)&  (B2)&  (B4)&   (B2)&  (B4)\notag \\
$15$         & (B4)         & (B2)&  (B4)&  (B2)&   (B4)&  (B2)\notag \\
$19$         & (B2)         & (B3)&  (B2)&  (B4)&   (B1)&  (B4)\notag \\
$23$         & (B4)         & (B2)&  (B4)&  (B2)&   (B4)&  (B2)\notag 
\end{tabular}
\]   
\caption{Cases (B1)--(B4) for $p$, $l$ taken modulo $24$} 
\label{Table33424}
\end{table}

The connection to $t_{p,l}$ is not as nice as in Section~\ref{CaseOne}.
We get the following values for $t_{p,l}$, if $p, l \equiv 3 \pmod{4}$
are distinct prime numbers less than $1000$.
\[
t_{p,l} \in
\begin{cases}
(\{ 4, 6, \ldots, 104 \}  
\cup \{ 110, 114, 122, 124, 132 \}) \setminus \{ 84, 88\}\, , 
& \text{in case (B1)}\\
\{ 0, 2, \ldots, 78\} \cup \{ 84, 100, 110 \} \, , & \text{in case (B2)}\\
\{ 0 \} \, , & \text{in case (B3)}\\
\{ 0 \} \, , & \text{in case (B4)} \, .
\end{cases}
\]
In general, i.e.\ without the restriction $p, l < 1000$, it is easy to see that
$t_{p,l}$ is always even. Moreover, it follows from \cite[Section~5]{Rattaggi3}
that $t_{p,l} = 0$ in the cases (B3), (B4), and $t_{p,l} > 0$ in case (B1).
The computations of $t_{p,l}$ combined with Conjecture~\ref{Conj334} lead to the following conjecture:
\begin{Conjecture} \label{Conjtpl}
Let $p,l \equiv 3 \pmod{4}$.
\begin{itemize}
\item[(1)] If $t_{p,l} = 0$, then $\Gamma_{p,l}^{ab} \cong \mathbb{Z}_2 \times \mathbb{Z}_8^2$  or 
$\Gamma_{p,l}^{ab} \cong \mathbb{Z}_2 \times \mathbb{Z}_3 \times \mathbb{Z}_4^2$ or 
$\Gamma_{p,l}^{ab} \cong \mathbb{Z}_2 \times \mathbb{Z}_4^2$.
\item[(2)] If $t_{p,l} = 2$, then $\Gamma_{p,l}^{ab} \cong \mathbb{Z}_2 \times \mathbb{Z}_8^2$.
\item[(3)] If $t_{p,l} \geq 4$, then $\Gamma_{p,l}^{ab} \cong \mathbb{Z}_2 \times \mathbb{Z}_3 \times \mathbb{Z}_8^2$
or $\Gamma_{p,l}^{ab} \cong \mathbb{Z}_2 \times \mathbb{Z}_8^2$.
\item[(4)] If $\Gamma_{p,l}^{ab} \cong \mathbb{Z}_2 \times \mathbb{Z}_3 \times \mathbb{Z}_4^2$ 
or $\Gamma_{p,l}^{ab} \cong \mathbb{Z}_2 \times \mathbb{Z}_4^2$, then $t_{p,l} = 0$.
\item[(5)] If $\Gamma_{p,l}^{ab} \cong \mathbb{Z}_2 \times \mathbb{Z}_3 \times \mathbb{Z}_8^2$,
then $t_{p,l} \geq 4$.
\end{itemize}
\end{Conjecture}

\section{The case $p \equiv 3 \pmod{4}$, $l \equiv 1 \pmod{4}$} \label{CaseThree}
The remaining case is $p \pmod{4} \ne l \pmod{4}$. Since $\Gamma_{p,l} \cong \Gamma_{l,p}$, we can
restrict to $p \equiv 3 \pmod{4}$, $l \equiv 1 \pmod{4}$.

\begin{Conjecture} \label{Conj314}
Let $p \equiv 3 \pmod{4}$, $l \equiv 1 \pmod{4}$.\\
If $l \equiv 1 \pmod{8}$, then
\[
\Gamma_{p,l}^{ab} \cong
\begin{cases}
\mathbb{Z}_2 \times \mathbb{Z}_3 \times \mathbb{Z}_8^2 \, ,  & 
\text{if }  \, p, l \equiv 1 \!\!\! \pmod{3} \; \text{\phantom{iii}(=: case (C1))} \\
\mathbb{Z}_2 \times \mathbb{Z}_8^2 \, ,                      & 
\text{else} \; \text{\phantom{iiiiiiiiiiiiiiiiiiiiiiiii}(=: case (C2))} \, .
\end{cases}
\]
If $l \equiv 5 \pmod{8}$, then
\[
\Gamma_{p,l}^{ab} \cong
\begin{cases}
\mathbb{Z}_2 \times \mathbb{Z}_3 \times \mathbb{Z}_4^2 \, ,  & 
\text{if }  \, p, l \equiv 1 \!\!\! \pmod{3} \; \text{\phantom{iii}(=: case (C3))} \\
\mathbb{Z}_2 \times \mathbb{Z}_4^2 \, ,                      & 
\text{else} \; \text{\phantom{iiiiiiiiiiiiiiiiiiiiiiiii}(=: case (C4))} \, .
\end{cases}
\]
\end{Conjecture}
Observe that the four conjectured possibilities for $\Gamma_{p,l}^{ab}$ are exactly the same
as in Conjecture~\ref{Conj334}.

See Table~\ref{Table31424} for the cases (C1)--(C4) expressed by $p$ and $l$ taken modulo $24$.
\begin{table}[ht]
\[
\begin{tabular}{r|rccccl}
             & $l \equiv 1$ & $5$ & $9$  & $13$ & $17$ & $21 \!\!\! \pmod{24}$ \notag \\ \hline
$p \equiv 3$ & (C2)         & (C4)&  (C2)&  (C4)&   (C2)&  (C4)\notag \\ 
$7$          & (C1)         & (C4)&  (C2)&  (C3)&   (C2)&  (C4)\notag \\
$11$         & (C2)         & (C4)&  (C2)&  (C4)&   (C2)&  (C4)\notag \\
$15$         & (C2)         & (C4)&  (C2)&  (C4)&   (C2)&  (C4)\notag \\
$19$         & (C1)         & (C4)&  (C2)&  (C3)&   (C2)&  (C4)\notag \\
$23$         & (C2)         & (C4)&  (C2)&  (C4)&   (C2)&  (C4)\notag 
\end{tabular}
\]   
\caption{Cases (C1)--(C4) for $p$, $l$ taken modulo $24$} 
\label{Table31424}
\end{table}

The behaviour of $t_{p,l}$ seems to be very similar as in Section~\ref{CaseTwo}.
We get the following values for $t_{p,l}$, if $p \equiv 3 \pmod{4}$, $l \equiv 1 \pmod{4}$
are prime numbers less than $1000$.
\[
t_{p,l} \in
\begin{cases}
(\{ 4, 6, \ldots, 48 \}  
\cup \{ 58 \}) \setminus \{ 40\}\, , 
& \text{in case (C1)}\\
\{ 0, 2, \ldots, 54\} \cup \{ 60 \} \, , & \text{in case (C2)}\\
\{ 0 \} \, , & \text{in case (C3)}\\
\{ 0 \} \, , & \text{in case (C4)} \, .
\end{cases}
\]
\begin{Conjecture}
Conjecture~\ref{Conjtpl} also holds if $p \equiv 3 \pmod{4}$, $l \equiv 1 \pmod{4}$.
\end{Conjecture}

\section{More conjectures} \label{MoreConj}
In this section, we give conjectures for the abelianization of the commutator subgroup
$[\Gamma_{p,l}, \Gamma_{p,l}]$ of $\Gamma_{p,l}$ and for a certain subgroup $\Lambda_{p,l}$ of
$\Gamma_{p,l}$ of index $4$ defined below.

\begin{Conjecture} \label{ConjComm114}
Let $p, l \equiv 1 \pmod{4}$.\\
If $p, l \equiv 1 \pmod{8}$, then
\[
[\Gamma_{p,l}, \Gamma_{p,l}]^{ab} \cong
\begin{cases}
\mathbb{Z}_2^2 \times \mathbb{Z}_{16}^2 \times \mathbb{Z}_{64} \, , & \text{if }  \, p, l \equiv 1 \!\!\! \pmod{3} \\
\mathbb{Z}_3 \times \mathbb{Z}_{16}^2 \times \mathbb{Z}_{64}  \, ,      & \text{else}\, .
\end{cases}
\]
If $p \equiv 5 \pmod{8}$ or $l \equiv 5 \pmod{8}$, then
\[
[\Gamma_{p,l}, \Gamma_{p,l}]^{ab} \cong
\begin{cases}
\mathbb{Z}_2^2 \times \mathbb{Z}_{16}^3 \, , & \text{if }  \, p, l \equiv 1 \!\!\! \pmod{3} \\
\mathbb{Z}_3 \times \mathbb{Z}_{16}^3 \, ,      & \text{else}\, .
\end{cases}
\]
\end{Conjecture}

\begin{Conjecture} \label{ConjComm334}
Let $p,l \equiv 3 \pmod{4}$.\\
If $p \! \pmod{8} \, = \, l  \! \pmod{8}$, then
\[
[\Gamma_{p,l}, \Gamma_{p,l}]^{ab} \cong
\begin{cases}
\mathbb{Z}_2^2 \times \mathbb{Z}_8^2 \times \mathbb{Z}_{64} \, ,  & \text{if } \, p, l \equiv 1 \!\!\! \pmod{3} \\
\mathbb{Z}_8^2 \times \mathbb{Z}_{64} \, , & \text{if } \, p = 3 \, \text{ or } \, l = 3\\
\mathbb{Z}_3 \times \mathbb{Z}_8^2 \times \mathbb{Z}_{64} \, ,       & \text{else}\, .
\end{cases}
\]
If $p  \! \pmod{8} \, \ne \, l  \! \pmod{8}$, then
\[
[\Gamma_{p,l}, \Gamma_{p,l}]^{ab} \cong
\begin{cases}
\mathbb{Z}_2^2 \times \mathbb{Z}_8^2 \times \mathbb{Z}_{16} \, , & \text{if }  \, p, l \equiv 1 \!\!\! \pmod{3} \\
\mathbb{Z}_8^2 \times \mathbb{Z}_{16}  & \text{if } \, p = 3 \, \text{ or } \, l = 3\\
\mathbb{Z}_3 \times \mathbb{Z}_8^2 \times \mathbb{Z}_{16} \, ,   & \text{else}\, .
\end{cases}
\]
\end{Conjecture}

The groups appearing in Conjecture~\ref{ConjComm314}
are again the same as in Conjecture~\ref{ConjComm334}:
\begin{Conjecture} \label{ConjComm314}
Let $p \equiv 3 \pmod{4}$ and $l \equiv 1 \pmod{4}$.\\
If $l \equiv 1 \pmod{8}$, then
\[
[\Gamma_{p,l}, \Gamma_{p,l}]^{ab} \cong
\begin{cases}
\mathbb{Z}_2^2 \times \mathbb{Z}_8^2 \times \mathbb{Z}_{64} \, , & \text{if } \, p, l \equiv 1 \!\!\! \pmod{3} \\
\mathbb{Z}_8^2 \times \mathbb{Z}_{64} \, , & \text{if } \, p = 3 \\
\mathbb{Z}_3 \times \mathbb{Z}_8^2 \times \mathbb{Z}_{64} \, ,   & \text{else}\, .
\end{cases}
\]
If $l \equiv 5 \pmod{8}$, then
\[
[\Gamma_{p,l}, \Gamma_{p,l}]^{ab} \cong
\begin{cases}
\mathbb{Z}_2^2 \times \mathbb{Z}_8^2 \times \mathbb{Z}_{16} \, , & \text{if }  \, p, l \equiv 1 \!\!\! \pmod{3} \\
\mathbb{Z}_8^2 \times \mathbb{Z}_{16} \, , & \text{if } \, p = 3 \\
\mathbb{Z}_3 \times \mathbb{Z}_8^2 \times \mathbb{Z}_{16} \, ,   & \text{else}\, .
\end{cases}
\]
\end{Conjecture}

Let $\Lambda_{p,l}$ be the following subgroup of $\Gamma_{p,l}$.
\[
\Lambda_{p,l} := \psi_{p,l}(\{x = x_0  + x_1 i + x_2 j + x_3 k \in \mathbb{H}(\mathbb{Z}); \,
x_0 \text{ odd}; \,
|x|^2 = p^{s} l^{t}, s,t \in 2\mathbb{N}_0 \}) \, .
\]
Observe that $\Lambda_{p,l}$ is the kernel of the surjective homomorphism 
$\Gamma_{p,l} \to \mathbb{Z}_2 \times \mathbb{Z}_2$ determined by
\[
\psi_{p,l}(x) \mapsto 
\begin{cases}
(1+2\mathbb{Z}, \, 0+2\mathbb{Z}) \, , & \text{if }  \, |x|^2 = p\\
(0+2\mathbb{Z}, \, 1+2\mathbb{Z}) \, , & \text{if }  \, |x|^2 = l \, ,
\end{cases}
\]
in particular $\Lambda_{p,l}$ is a normal subgroup of $\Gamma_{p,l}$ of index $4$.
It seems that the abelianization of $\Lambda_{p,l}$ does not depend on $p$ and $l$, if $p,l > 3$.
\begin{Conjecture} \label{ConjectureLambda}
Let $p,l$ be any pair of distinct odd prime numbers. Then
\[
\Lambda_{p,l}^{ab} \cong
\begin{cases}
\mathbb{Z}_2 \times \mathbb{Z}_8^2 \, ,       & \text{if }  \, p = 3 \, \text{ or } \, l = 3\\
\mathbb{Z}_2 \times \mathbb{Z}_3 \times \mathbb{Z}_8^2 \, ,   & \text{else}\, .
\end{cases}
\]
\end{Conjecture}

Since the conjectured abelianizations of the groups $\Gamma_{p,l}$, $[\Gamma_{p,l}, \Gamma_{p,l}]$
and $\Lambda_{p,l}$ are never $2$-generated, we also get the following conjecture:
\begin{Conjecture}
Let $p,l$ be any pair of distinct odd prime numbers. Then the groups
$\Gamma_{p,l}$, $[\Gamma_{p,l}, \Gamma_{p,l}]$
and $\Lambda_{p,l}$ are not $2$-generated.
\end{Conjecture}


\begin{thebibliography}{99}
\bibitem{KR}
  Kimberley, Jason S.; Robertson, Guyan,
  \emph{Groups acting on products of trees, tiling systems and analytic K-theory},
  New York J. Math. \textbf{8}(2002), 111--131 (electronic).
\bibitem{Rattaggi}
  Rattaggi, Diego,
  \emph{Computations in groups acting on a product of trees: 
   normal subgroup structures and quaternion lattices},
   Ph.D.~thesis, ETH Z\"urich, 2004.
\bibitem{Rattaggi2}
  Rattaggi, Diego,
  \emph{Anti-tori in square complex groups},
  to appear in Geom.\ Dedicata.
  A preprint is available at arXiv:math.GR/0411547.
\bibitem{Rattaggi3}
  Rattaggi, Diego,
  \emph{On infinite groups generated by two quaternions},
  Preprint 2005, see arXiv:math.GR/0502512.
\bibitem{RaRo}
  Rattaggi, Diego; Robertson, Guyan,
  \emph{Abelian subgroup structure of square complex groups 
  and arithmetic of quaternions},
  J. Algebra  \textbf{286}(2005),  no.~1, 57--68.
\end{thebibliography}
\end{document}